\documentclass[11pt,twoside]{article}
\usepackage{latexsym,amsmath,amsthm,amssymb,amsfonts}
\usepackage{mathrsfs}
\usepackage{amscd}
\usepackage{cite}

\usepackage[ansinew]{inputenc}
\numberwithin{equation}{section}

\def\proof{{\bf Proof.}}

\catcode`@=11

\newskip\plaincentering \plaincentering=0pt plus 1000pt minus 1000pt
\def\@plainlign{\tabskip=0pt\everycr={}}
\def\eqalignno#1{\displ@y \tabskip\plaincentering
  \halign to\displaywidth{\hfil$\@lign\displaystyle{##}$\tabskip\z@skip
    &$\@lign\displaystyle{{}##}$\hfil\tabskip\plaincentering
    &\llap{$\@lign##$}\tabskip\z@skip\crcr
    #1\crcr}}
\def\leqalignno#1{\displ@y \tabskip\plaincentering
  \halign to\displaywidth{\hfil$\@lign\displaystyle{##}$\tabskip\z@skip
    &$\@lign\displaystyle{{}##}$\hfil\tabskip\plaincentering
    &\kern-\displaywidth\rlap{$\@lign##$}\tabskip\displaywidth\crcr
    #1\crcr}}
\def\plainLet@{\relax\iffalse{\fi\let\\=\cr\iffalse}\fi}
\def\plainvspace@{\def\vspace##1{\noalign{\vskip##1}}}

\def\intic@{\mathchoice{\hskip5\p@}{\hskip4\p@}{\hskip4\p@}{\hskip4\p@}}
\def\negintic@
 {\mathchoice{\hskip-5\p@}{\hskip-4\p@}{\hskip-4\p@}{\hskip-4\p@}}
\def\intkern@{\mathchoice{\!\!\!}{\!\!}{\!\!}{\!\!}}
\def\intdots@{\mathchoice{\cdots}{{\cdotp}\mkern1.5mu
    {\cdotp}\mkern1.5mu{\cdotp}}{{\cdotp}\mkern1mu{\cdotp}\mkern1mu
      {\cdotp}}{{\cdotp}\mkern1mu{\cdotp}\mkern1mu{\cdotp}}}
\newcount\intno@
\def\iint{\intno@=\tw@\futurelet\next\ints@}
\def\iiint{\intno@=\thr@@\futurelet\next\ints@}
\def\iiiint{\intno@=4 \futurelet\next\ints@}
\def\idotsint{\intno@=\z@\futurelet\next\ints@}
\def\ints@{\findlimits@\ints@@}
\newif\iflimtoken@
\newif\iflimits@
\def\findlimits@{\limtoken@false\limits@false\ifx\next\limits
 \limtoken@true\limits@true\else\ifx\next\nolimits\limtoken@true\limits@false
    \fi\fi}
\def\multintlimits@{\intop\ifnum\intno@=\z@\intdots@
  \else\intkern@\fi
    \ifnum\intno@>\tw@\intop\intkern@\fi
     \ifnum\intno@>\thr@@\intop\intkern@\fi\intop}
\def\multint@{\int\ifnum\intno@=\z@\intdots@\else\intkern@\fi
   \ifnum\intno@>\tw@\int\intkern@\fi
    \ifnum\intno@>\thr@@\int\intkern@\fi\int}
\def\ints@@{\iflimtoken@\def\ints@@@{\iflimits@
   \negintic@\mathop{\intic@\multintlimits@}\limits\else
    \multint@\nolimits\fi\eat@}\else
     \def\ints@@@{\multint@\nolimits}\fi\ints@@@}
\def\Sb{_\bgroup\vspace@
        \baselineskip=\fontdimen10 \scriptfont\tw@
        \advance\baselineskip by \fontdimen12 \scriptfont\tw@
        \lineskip=\thr@@\fontdimen8 \scriptfont\thr@@
        \lineskiplimit=\thr@@\fontdimen8 \scriptfont\thr@@
        \Let@\vbox\bgroup\halign\bgroup \hfil$\scriptstyle
            {##}$\hfil\cr}
\def\endSb{\crcr\egroup\egroup\egroup}
\def\Sp{^\bgroup\vspace@
        \baselineskip=\fontdimen10 \scriptfont\tw@
        \advance\baselineskip by \fontdimen12 \scriptfont\tw@
        \lineskip=\thr@@\fontdimen8 \scriptfont\thr@@
        \lineskiplimit=\thr@@\fontdimen8 \scriptfont\thr@@
        \Let@\vbox\bgroup\halign\bgroup \hfil$\scriptstyle
            {##}$\hfil\cr}
\def\endSp{\crcr\egroup\egroup\egroup}
\def\Let@{\relax\iffalse{\fi\let\\=\cr\iffalse}\fi}
\def\vspace@{\def\vspace##1{\noalign{\vskip##1 }}}
\def\aligned{\,\vcenter\bgroup\plainvspace@\plainLet@\openup\jot\m@th\ialign
  \bgroup \strut\hfil$\displaystyle{##}$&$\displaystyle{{}##}$\hfil\crcr}
\def\endaligned{\crcr\egroup\egroup}
\def\matrix{\,\vcenter\bgroup\plainLet@\plainvspace@
    \normalbaselines
  \m@th\ialign\bgroup\hfil$##$\hfil&&\quad\hfil$##$\hfil\crcr
    \mathstrut\crcr\noalign{\kern-\baselineskip}}
\def\endmatrix{\crcr\mathstrut\crcr\noalign{\kern-\baselineskip}\egroup
                \egroup\,}
\newtoks\hashtoks@
\hashtoks@={#}
\def\format{\crcr\egroup\iffalse{\fi\ifnum`}=0 \fi\format@}
\def\format@#1\\{\def\preamble@{#1}%
  \def\c{\hfil$\the\hashtoks@$\hfil}%
  \def\r{\hfil$\the\hashtoks@$}%
  \def\l{$\the\hashtoks@$\hfil}%
  \setbox\z@=\hbox{\xdef\Preamble@{\preamble@}}\ifnum`{=0 \fi\iffalse}\fi
   \ialign\bgroup\span\Preamble@\crcr}

\def\cases{\left\{\,\vcenter\bgroup\plainvspace@
     \normalbaselines\openup\jot\m@th
      \plainLet@\ialign\bgroup$\displaystyle{##}$\hfil&\quad$\displaystyle{{}##}$\hfil\crcr
      \mathstrut\crcr\noalign{\kern-\baselineskip}}

\newif\iftagsleft@
\tagsleft@true
\def\TagsOnRight{\global\tagsleft@false}
\def\tag#1$${\iftagsleft@\leqno\else\eqno\fi
 \hbox{\def\pagebreak{\global\postdisplaypenalty-\@M}%
 \def\nopagebreak{\global\postdisplaypenalty\@M}\rm(#1\unskip)}%
  $$\postdisplaypenalty\z@\ignorespaces}
\interdisplaylinepenalty=\@M
\def\plainallowdisplaybreak@{\def\allowdisplaybreak{\noalign{\allowbreak}}}
\def\plaindisplaybreak@{\def\displaybreak{\noalign{\break}}}
\def\align#1\endalign{\def\tag{&}\plainvspace@\plainallowdisplaybreak@\plaindisplaybreak@
  \iftagsleft@\plainlalign@#1\endalign\else
   \plainralign@#1\endalign\fi}
\def\plainralign@#1\endalign{\displ@y\plainLet@\tabskip\plaincentering\halign to\displaywidth
     {\hfil$\displaystyle{##}$\tabskip=\z@&$\displaystyle{{}##}$\hfil
       \tabskip=\plaincentering&\llap{\hbox{\rm(##\unskip)}}\tabskip\z@\crcr
             #1\crcr}}
\def\plainlalign@
 #1\endalign{\displ@y\plainLet@\tabskip\plaincentering\halign to \displaywidth
   {\hfil$\displaystyle{##}$\tabskip=\z@&$\displaystyle{{}##}$\hfil
   \tabskip=\plaincentering&\kern-\displaywidth
        \rlap{\hbox{\rm(##\unskip)}}\tabskip=\displaywidth\crcr
               #1\crcr}}

\def\re@#1{\par\hangindent\parindent\indent\llap{#1\enspace}\ignorespaces}
\def\qfootnote#1{\edef\@sf{\spacefactor\the\spacefactor}{}#1\@sf
      \insert\footins{\let\egroup=}\footnotesize 
      \interlinepenalty100 \let\par=\endgraf
        \leftskip=0pt \rightskip=0pt
        \splittopskip=10pt plus 1pt minus 1pt \floatingpenalty=20000
   \smallskip\re@{#1}\bgroup\strut\aftergroup{\strut\egroup}\let\next}
\catcode`\@=\active
\TagsOnRight
\def\proof{{\bf Proof.}\ }

\begin{document}
\title{\bf Estimates on the first two buckling eigenvalues on spherical domains
\footnote{This research is supported by National Natural Science
Foundation of China (Project No. 10671181) and National Natural
Science Foundation of Henan (Project No. 092300410143;
2009A110010).}}
\author{Guangyue Huang,\ Xingxiao Li
\footnote{The corresponding author. Email: xxl$@$henannu.edu.cn},\ Xuerong Qi\\
{\normalsize Department of Mathematics, Henan Normal University}
\\{\normalsize Xinxiang 453007, Henan, P.R. China}}
\date{}
\maketitle
\begin{quotation}
\noindent{\bf Abstract.}~In this paper, we study the first two
eigenvalues of the buckling problem on spherical domains. We obtain
an estimate on the second eigenvalue in terms of the first
eigenvalue, which improves one recent result
obtained by Wang-Xia in \cite{wangxia07}.

\end{quotation}

\section{Introduction}

Let $\Omega$ be a connected bounded domain in an $n$-dimensional
Euclidean space $\mathbb{R}^n$ and $\nu$ be the unit outward normal
vector field of $\partial\Omega$. The well-known eigenvalue problem
\begin{equation}\label{Intro1}\left\{\begin{array}{ll}
\Delta^2 u=-\Lambda\Delta u & {\rm in}\ \Omega,\\
u={\partial u\over\partial\nu}=0 & {\rm on}\
\partial\Omega
\end{array}\right.
\end{equation}
is called a {\em buckling problem}, which is used to describe the
critical buckling load of a clamped plate subjected to a uniform
compressive force around its boundary.

Let $$0<\Lambda_1\leq\Lambda_2\leq \Lambda_3\leq \cdots$$ denote the
successive eigenvalues for \eqref{Intro1}, where each eigenvalue is
repeated according to its multiplicity. In 1956,
Payne-P\'{o}lya-Weinberger \cite{ppw2} proved that
$$
\Lambda_2\leq \left(1+{4\over n}\right)\Lambda_1.
$$
Subsequently, Hile-Yeh \cite{hile1} improved the above inequality as
follows:
$$
\Lambda_2\leq {n^2+8n+20\over (n+2)^2}\Lambda_1.
$$ In 1994, Chen-Qian \cite{chenqian94} considered the following more general
eigenvalue problem: $$\left\{\begin{array}{ll}
(-\Delta)^pu=\Lambda(-\Delta)^qu & {\rm in}\ \Omega,\\
u={\partial u\over\partial\nu}=\cdots={\partial^{p-1}
u\over\partial\nu^{p-1}}=0 & {\rm on}\
\partial\Omega
\end{array}\right.
$$
with $p$ and $q$ are positive integers and $p>q$. They proved that
$$
\Lambda_2\leq {(n+2q)^2+4p(2p+n-2)-4q(2q+n-2)\over
(n+2q)^2}\Lambda_1.
$$
In 1998, Ashbaugh \cite{ash} found that
$$\sum_{i=1}^n\Lambda_{i+1}\leq (n+4)\Lambda_1.$$
For answering a question of Ashbaugh in \cite{ash}, Cheng-Yang
\cite{chengyang06} obtained in a recent survey paper a universal
inequality for higher eigenvalues of the buckling problem. In fact,
they proved that
$$\sum_{i=1}^k(\Lambda_{k+1}-\Lambda_i)^2\leq{4(n+2)
\over n^2}\sum_{i=1}^k(\Lambda_{k+1}-\Lambda_i)\Lambda_i.$$

In 2007, Wang-Xia \cite{wangxia07} considered the buckling problem
on domains in a unit sphere and obtained the following result:

{\bf Theorem A.} {\em Let $\Lambda_i$ be the $i^{th}$ eigenvalue of
the following eigenvalue problem:
$$\Delta^2 u=-\Lambda\Delta u \ \ {\rm in}\ \Omega,\ \ \ \
u={\partial u\over\partial\nu}=0\ \ {\rm on}\
\partial\Omega,$$ where $\Omega$ is a connected domain in an $n$-dimensional unit sphere with
smooth boundary $\partial\Omega$ and $\nu$ is the unit outward
normal vector field of $\partial\Omega$. Then for any $\delta>0$, it
holds that
\begin{equation}\label{Intro2}
\aligned 2\sum_{i=1}^k(\Lambda_{k+1}-\Lambda_i)^2\leq&
\sum_{i=1}^k(\Lambda_{k+1}-\Lambda_i)^2\left(\delta\Lambda_i
+{\delta^2(\Lambda_i-(n-2))\over4(\delta\Lambda_i+n-2)}\right)\\
&+
{1\over\delta}\sum_{i=1}^k(\Lambda_{k+1}-\Lambda_i)\left(\Lambda_i+{(n-2)^2\over4}\right).
\endaligned
\end{equation}
}

Recently, Huang-Li-Cao \cite{huang} improved the above result as
follows:

{\bf Theorem B.} {\em Under the assumption of Theorem A. Then for
any $\delta>0$,
\begin{equation}\label{Intro3}
\aligned &\sum_{i=1}^k(\Lambda_{k+1}-\Lambda_i)^2
\left(2+\frac{n-2}{\Lambda_i-(n-2)}\right)\\
\leq& 2\left\{\sum_{i=1}^k(\Lambda_{k+1}-\Lambda_i)^2\left(\Lambda_i
-\frac{n-2}{\Lambda_i-(n-2)}\right)\right\}^{\frac{1}{2}}\\
&\qquad \qquad\qquad
\times\left\{\sum_{i=1}^k(\Lambda_{k+1}-\Lambda_i)
\left(\Lambda_i+\frac{(n-2)^2}{4}\right)\right\}^{\frac{1}{2}}.
\endaligned
\end{equation}
}

In the present paper, we consider the first two eigenvalues of the
buckling problem on spherical domains and obtain the following
result:

{\bf Theorem 1.1.} {\em Let $\Lambda_i$ be the $i^{th}$ eigenvalue
of the following eigenvalue problem:
\begin{equation}\label{Intro4}\Delta^2 u=-\Lambda\Delta u \ \ {\rm in}\ \Omega,\ \ \ \
u={\partial u\over\partial\nu}=0\ \ {\rm on}\
\partial\Omega,\end{equation} where $\Omega$ is a connected domain
in an $n$-dimensional unit sphere with smooth boundary
$\partial\Omega$ and $\nu$ is the unit outward normal vector field
of $\partial\Omega$. Then we have
\begin{equation}\label{Intro5}\aligned
\Lambda_2\leq\Lambda_1+\left({ n(n-\Lambda_1)\over\Lambda_1}
+2(n+2)\right){4\Lambda_1+(n-2)^2\over(n+2)^2}.
\endaligned
\end{equation}
}

{\bf Corollary 1.2.} {\em Under the same assumption of Theorem 1.1,
we have
\begin{equation}\label{Intro6}
\Lambda_2\leq\left(1+{8\over n+2}\right)\Lambda_1+{2(n-2)^2\over
n+2}.
\end{equation}
}

{\bf Remark.} From the inequality \eqref{2sec18} in section 2, we
have $$\Lambda_1\geq n.$$ Hence, we derive from \eqref{Intro3} that
\begin{equation}\label{Intro7}
(\Lambda_2-\Lambda_1)\leq\Lambda_1
^{\frac{1}{2}}\left\{(\Lambda_{2}-\Lambda_1)
\left(\Lambda_1+\frac{(n-2)^2}{4}\right)\right\}^{\frac{1}{2}}.
\end{equation} That is,
\begin{equation}\label{Intro8}
\Lambda_2\leq\Lambda_1+\Lambda_1\left(\Lambda_1+{(n-2)^2\over4}\right).
\end{equation}
Note that $$\left(1+{8\over n+2}\right)\Lambda_1+{2(n-2)^2\over
n+2}<\Lambda_1+\Lambda_1\left(\Lambda_1+{(n-2)^2\over4}\right).$$
Thus it is not hard to see that, for the first two eigenvalues of
problem \eqref{Intro4}, the inequality \eqref{Intro6} is sharper
than \eqref{Intro3}.

\section{Proof of Theorem}

Let $x_1,x_2,\ldots, x_{n+1}$ be the standard coordinate functions
of the Euclidean space $\mathbb{R}^{n+1}$. Define
$$\mathbb{S}^n=\left\{ (x_1,x_2,\ldots, x_{n+1})\in\mathbb{R}^{n+1}\ ;\
\sum_{\alpha=1}^{n+1}x_\alpha^2=1\right\}.$$ Denote by $u_1$ the
eigenfunction corresponding to $\Lambda_1$ of the eigenvalue problem
\eqref{Intro4}, that is
\begin{equation*}\left\{\begin{array}{l}
\Delta^2 u_1=-\Lambda_1\Delta u_1 \ \ {\rm in}\ \Omega,\\
u_1={\partial u_1\over\partial\nu}=0\ \ {\rm on}\
\partial\Omega,\\
\int_\Omega  \langle\nabla u_1,\nabla u_1\rangle=1.
\end{array}\right.
\end{equation*}
Let
\begin{equation}\label{2sec1} \varphi_{\alpha}=x_\alpha u_1-C_\alpha u_1,
\end{equation}
where $$C_\alpha=\int_\Omega x_\alpha u_1(-\Delta)u_1.$$  Then one
gets $\varphi_{\alpha}=\partial\varphi_{\alpha}/\partial \nu=0$ on
$\partial\Omega$ and
\begin{equation}\label{2sec2}\int_\Omega \langle\nabla\varphi_{\alpha},\nabla u_1\rangle=0.\end{equation}
It follows from the Rayleigh-Ritz inequality that
\begin{equation}\label{2sec3}\Lambda_2\leq
{\int_\Omega \varphi_{\alpha}\Delta^2\varphi_{\alpha}\over
\int_\Omega |\nabla\varphi_{\alpha}|^2}.\end{equation}

Using integration by parts and the definition of $\varphi_\alpha$,
one finds that
\begin{equation}\aligned
\label{2sec4}\int_\Omega |\nabla\varphi_{\alpha}|^2=&\int_\Omega
\varphi_{\alpha}(-\Delta)\varphi_{\alpha}\\
=&\int_\Omega \varphi_{\alpha}(-\Delta)(x_\alpha u_1)\\
=&\int_\Omega \varphi_{\alpha}x_\alpha(-\Delta) u_1-\int_\Omega
\varphi_{\alpha}[\Delta(x_\alpha u_1)+x_\alpha(-\Delta) u_1].
\endaligned
\end{equation}
Note that
\begin{equation}\label{2sec5}\aligned
\int_\Omega \varphi_{\alpha}\Delta^2\varphi_{\alpha}
=&\int_\Omega \varphi_{\alpha}\Delta^2(x_\alpha u_1)\\
=&\Lambda_1\int_\Omega \varphi_{\alpha}x_\alpha(-\Delta) u_1
+\int_\Omega \varphi_{\alpha}[\Delta^2(x_\alpha
u_1)-x_\alpha\Delta^2 u_1].
\endaligned
\end{equation}
It follows from \eqref{2sec3}, \eqref{2sec4} and \eqref{2sec5} that
\begin{equation}\label{2sec6}\aligned
(\Lambda_2-\Lambda_1)\int_\Omega |\nabla\varphi_{\alpha
}|^2\leq&\Lambda_1\int_\Omega\varphi_{\alpha}[\Delta(x_\alpha u_1)+x_\alpha(-\Delta)u_1]\\
&+\int_\Omega\varphi_{\alpha}[\Delta^2(x_\alpha
u_1)-x_\alpha\Delta^2 u_1].
\endaligned
\end{equation}
Again, by using integration by parts, one gets that
\begin{equation}\label{2sec7}\aligned
\Lambda_1&\sum_{\alpha=1}^{n+1}\int_\Omega\varphi_{\alpha}[\Delta(x_\alpha
u_1)+x_\alpha(-\Delta)u_1]\\
=&\Lambda_1\sum_{\alpha=1}^{n+1}\int_\Omega x_{\alpha
}u_1[\Delta(x_\alpha u_1)+x_\alpha(-\Delta) u_1]\\
=&\Lambda_1\sum_{\alpha=1}^{n+1}\int_\Omega x_{\alpha
}u_1[-nx_{\alpha
}u_1+2\langle\nabla x_\alpha,\nabla u_1\rangle]\\
=&-n\Lambda_1\int_\Omega u_1^2.
\endaligned\end{equation} Inserting \eqref{2sec7} into \eqref{2sec6}
yields
\begin{equation}\label{2sec8}
(\Lambda_2-\Lambda_1)\int_\Omega |\nabla\varphi_{\alpha
}|^2\leq-n\Lambda_1\int_\Omega u_1^2\\
+\int_\Omega\varphi_{\alpha}[\Delta^2(x_\alpha u_1)-x_\alpha\Delta^2
u_1].
\end{equation}

Let $\nabla$ and $\nabla^2$ be the gradient operator and the Hessian
operator of $\mathbb{S}^n$, respectively.  Then we have
\begin{equation}\label{2sec9}
\Delta(x_\alpha  u_1)=-nx_\alpha u_1+2\langle\nabla x_\alpha,\nabla
u_1\rangle+x_\alpha\Delta u_1,
\end{equation} where we have used $\Delta x_\alpha=-n x_\alpha$. From the Ricci identity,
it follows that \begin{equation}\label{2sec10}\Delta\langle\nabla
x_\alpha,\nabla u_1\rangle=2\langle\nabla^2 x_\alpha,\nabla^2
u_1\rangle+\langle\nabla x_\alpha,\nabla \Delta
u_1\rangle+(n-2)\langle\nabla x_\alpha,\nabla u_1\rangle.
\end{equation}
It is well-known that $$\nabla^2 x_\alpha=-x_\alpha\langle\ , \
\rangle.$$ Therefore, \eqref{2sec10} becomes
\begin{equation}\label{2sec11}
\Delta\langle\nabla x_\alpha,\nabla u_1\rangle=-2x_\alpha\Delta
u_1+\langle\nabla x_\alpha,\nabla \Delta
u_1\rangle+(n-2)\langle\nabla x_\alpha,\nabla u_1\rangle.
\end{equation}
By virtue of \eqref{2sec9} and \eqref{2sec11}, a direct calculation
yields
\begin{equation}\label{2sec12}
 \Delta^2(x_\alpha u_1)= x_\alpha \Delta^2
u_1+n^2x_\alpha u_1-2(n+2)x_\alpha\Delta u_1-4\langle\nabla
x_\alpha,\nabla u_1\rangle+4\langle\nabla x_\alpha,\nabla\Delta
u_1\rangle.
\end{equation}
Then
\begin{equation}\label{2sec13}\aligned
\sum_{\alpha=1}^{n+1}&\int_\Omega\varphi_{\alpha}[\Delta^2(x_\alpha
u_1)-x_\alpha\Delta^2
u_1]\\
=&\sum_{\alpha=1}^{n+1}\int_\Omega x_{\alpha}u_1[\Delta^2(x_\alpha
u_1)-x_\alpha\Delta^2 u_1]\\
=&\sum_{\alpha=1}^{n+1}\int_\Omega x_{\alpha}u_1[n^2x_\alpha
u_1-2(n+2)x_\alpha\Delta u_1-4\langle\nabla x_\alpha,\nabla
u_1\rangle+4\langle\nabla x_\alpha,\nabla\Delta u_1\rangle]\\
=&n^2\int_\Omega u_1^2-2(n+2)\int_\Omega u_1\Delta u_1\\
=&n^2\int_\Omega u_1^2+2(n+2).
\endaligned\end{equation}
From \eqref{2sec8}  and  \eqref{2sec13}, we arrive at
\begin{equation}\label{2sec14}
(\Lambda_2-\Lambda_1)\sum_{\alpha=1}^{n+1}\int_\Omega
|\nabla\varphi_{\alpha }|^2\leq n(n-\Lambda_1)\int_\Omega
u_1^2+2(n+2).
\end{equation}

Let
$$D_\alpha=\int_\Omega\langle\nabla\varphi_\alpha,\nabla\langle\nabla
x_\alpha,\nabla u_1\rangle-{n-2\over2}x_\alpha\nabla u_1\rangle.$$
Then $$\aligned
\sum_{\alpha=1}^{n+1}D_\alpha=&\sum_{\alpha=1}^{n+1}\int_\Omega\langle\nabla
(x_\alpha u_1),\nabla\langle\nabla x_\alpha,\nabla
u_1\rangle-{n-2\over2}x_\alpha\nabla
u_1\rangle\\
&-\sum_{\alpha=1}^{n+1}C_\alpha\int_\Omega\langle\nabla
u_1,\nabla\langle\nabla
x_\alpha,\nabla u_1\rangle-{n-2\over2}x_\alpha\nabla u_1\rangle\\
=&\sum_{\alpha=1}^{n+1}\int_\Omega\langle\nabla (x_\alpha
u_1),\nabla\langle\nabla x_\alpha,\nabla
u_1\rangle-{n-2\over2}x_\alpha\nabla
u_1\rangle\\
=&\sum_{\alpha=1}^{n+1}\int_\Omega\langle\nabla (x_\alpha
u_1),\nabla\langle\nabla x_\alpha,\nabla
u_1\rangle\rangle-{n-2\over2}\sum_{\alpha=1}^{n+1}\int_\Omega\langle\nabla
(x_\alpha u_1),x_\alpha\nabla
u_1\rangle\\
=&-2-{n-2\over2}\\
=&-{n+2\over2}.
\endaligned$$
Since $$\aligned
&\sum_{\alpha=1}^{n+1}\int_\Omega|\nabla\langle\nabla
x_\alpha,\nabla u_1\rangle-{n-2\over2}x_\alpha\nabla
u_1|^2\\
=&\sum_{\alpha=1}^{n+1}\left(\int_\Omega|\nabla\langle\nabla
x_\alpha,\nabla
u_1\rangle|^2-(n-2)\int_\Omega\langle\nabla\langle\nabla
x_\alpha\nabla u_1\rangle,x_\alpha\nabla
u_1\rangle\right.\\
&\left.+\frac{(n-2)^2}{4}\int_\Omega|x_\alpha\nabla
u_1|^2\right)\\
=&\Lambda_1+{(n-2)^2\over4}.\endaligned$$ It follows that
\begin{equation}\label{2sec15}\aligned
{(n+2)^2\over4}=&\left(\sum_{\alpha=1}^{n+1}D_\alpha\right)^2\\
\leq&\left(\sum_{\alpha=1}^{n+1}\int_\Omega|\nabla\varphi_\alpha|^2\right)
\left(\sum_{\alpha=1}^{n+1}\int_\Omega|\nabla\langle\nabla
x_\alpha,\nabla u_1\rangle-{n-2\over2}x_\alpha\nabla
u_1|^2\right)\\
=&\left(\Lambda_1+{(n-2)^2\over4}\right)
\sum_{\alpha=1}^{n+1}\int_\Omega|\nabla\varphi_\alpha|^2.
\endaligned
\end{equation}
By \eqref{2sec15}, we obtain
\begin{equation}\label{2sec16}
\sum_{\alpha=1}^{n+1}\int_\Omega|\nabla\varphi_\alpha|^2
\geq{(n+2)^2\over 4\Lambda_1+(n-2)^2}.\end{equation} Applying
\eqref{2sec16} to \eqref{2sec14} yields
\begin{equation}\label{2sec17}
\Lambda_2\leq\Lambda_1+\left(n(n-\Lambda_1)\int_\Omega
u_1^2+2(n+2)\right){4\Lambda_1+(n-2)^2\over(n+2)^2}.
\end{equation}

{\bf Lemma 2.1.} {\em Let $\Omega_1$, $\Omega_2$ be two connected
bounded domains in $\Bbb{S}^n$ and $\Omega_1\subset\Omega_2$. Then
it holds that $\Lambda_1(\Omega_1)\geq \Lambda_1(\Omega_2)$.}

\proof Let $u_1(\Omega_1)$ be the eigenfunction corresponding to
$\Lambda_1(\Omega_1)$. Then the function defined by $$\widetilde{
u_1}=\left\{\begin{array}{ll} u_1 &\ \ \  {\rm in }\ \ \Omega_1,\\
0 &\ \ \  {\rm in }\ \ \Omega_2- \Omega_1
\end{array}\right.
$$
is a eigenfunction corresponding to $\Lambda_1(\Omega_2)$. This
easily proves Lemma 2.1.

Let $\lambda_1$ be the first eigenvalue of Laplacian. Then
$\Lambda_1(\Bbb{S}^n)=\lambda_1(\Bbb{S}^n)=n$ because there are no
boundary conditions in this case. It follows that
\begin{equation}\label{2sec18}\Lambda_1(\Omega_1)\geq
n\end{equation} by setting $\Omega=\Omega_1$ and
$\Omega_2=\Bbb{S}^n$ in Lemma 2.1.

Using Schwarz inequality, we have
$$
1=\int_\Omega |\nabla u_1|^2=\int_\Omega
u_1(-\Delta)u_1\leq\left(\int_\Omega u_1^2\int_\Omega(\Delta
u_1)^2\right)^{1/2}=\left(\Lambda_1\int_\Omega u_1^2\right)^{1/2}.
$$
Hence,
\begin{equation}\label{2sec19}\int_\Omega u_1^2\geq {1\over \Lambda_1}.\end{equation}

Applying \eqref{2sec18} and \eqref{2sec19} into \eqref{2sec17}
yields
$$
\Lambda_2\leq\Lambda_1+\left({ n(n-\Lambda_1)\over\Lambda_1}
+2(n+2)\right){4\Lambda_1+(n-2)^2\over(n+2)^2}.
$$
This completes the proof of Theorem 1.1.

{\em Proof of Corollary 1.2.}\ From the inequality \eqref{2sec18},
we have
\begin{equation}\label{2sec20}{n(n-\Lambda_1)\over\Lambda_1}+2(n+2)\leq
2(n+2).\end{equation} Applying \eqref{2sec20} to \eqref{Intro5}
completes the proof of Corollary 1.2.

\end{document}
\end{document}